\documentclass[12pt]{article}%
\usepackage[intlimits]{amsmath}
\usepackage{amssymb}
\providecommand{\U}[1]{\protect\rule{.1in}{.1in}}
\providecommand{\U}[1]{\protect\rule{.1in}{.1in}}
\newtheorem{thm}{Theorem}[section]

\newtheorem{lem}[thm]{Lemma}

\newtheorem{prop}[thm]{Proposition}
\newtheorem{rem}[thm]{Remark}

\newenvironment{pf}[1][Proof]{\noindent\textbf{#1.} }{\hfill$\Box$\vspace{0.4cm}}

\begin{document}


\thispagestyle{empty}
\title{A probabilistic approach to the asymptotics of the length of the longest alternating subsequence.}
\author{Christian Houdr\'{e} \thanks{Georgia Institute of Technology, School of
Mathematics, Atlanta, Georgia, 30332, USA, houdre@math.gatech.edu. Supported
in part by the NSA grant H98230-09-1-0017.}
\and Ricardo Restrepo \thanks{Georgia Institute of Technology, School of
Mathematics, Atlanta, Georgia, 30332, USA, restrepo@math.gatech.edu.}
\thanks{Universidad de Antioquia, Departamento de Matematicas, Medellin,
Colombia.
}}

\maketitle


\begin{abstract}
\noindent Let $LA_{n}(\tau)$ be the length of the longest alternating
subsequence of a uniform random permutation $\tau\in\left[  n\right]  $. 
Classical probabilistic
arguments are used to rederive the asymptotic mean, variance and limiting
law of $LA_{n}\left(  \tau\right)  $. Our methodology is robust enough to
tackle similar problems for finite alphabet random words or even Markovian sequences in which case our results are mainly original. A sketch of how some cases of pattern restricted permutations can also be tackled with probabilistic methods is finally presented.

\end{abstract}



\section{Introduction}

Let $a:=(a_{1},a_{2},\ldots,a_{n})$ be a sequence of length $n$ whose elements
belong to a totally ordered set $\Lambda$. Given an increasing set of indices
$\{\ell_{i}\}_{i=1}^{m}$, we say that the subsequence $(a_{\ell_{1}}%
,a_{\ell_{2}},\ldots,a_{\ell_{m}})$ is alternating if $a_{\ell_{1}}%
>a_{\ell_{2}}<a_{\ell_{3}}>\cdots a_{\ell_{m}}$. The \emph{length of the
longest alternating subsequence} is then defined as%
\[
\operatorname*{LA}\nolimits_{n}(a):=\max\left\{  m:a\text{ has an alternating
subsequence of length }m\right\}  \text{.}%
\]

We revisit, here, the problem of finding the asymptotic behavior (in mean,
variance and limiting law) of the length of the longest alternating
subsequence in the context of random permutations and random words. For random
permutations, these problems have seen complete solutions with contributions
independently given (in alphabetical order) by Pemantle,
Stanley and Widom. The reader will find in \cite{stanley2006increasing}
a comprehensive survey, with precise
bibliography and credits, on these and related problems. 
In the context of random words, Mansour
\cite{mansour2008longest} contains very recent
contributions where mean and variance are obtained. Let us just say that, to date, the
proofs developed to solve these problems are of a combinatorial or analytic
nature and that we wish below to provide probabilistic ones. Our approach is
developed via iid sequences uniformly distributed on $[0,1]$, counting minima
and maxima and the central limit theorem for $2$-dependent random variables.
Not only does our approach recover the permutation case, but it works as well for random words, $a\in \mathcal{A}^n$ where $\mathcal{A}$ is a finite ordered alphabet, recovering known results and providing new ones. Properly modified it also works for several kinds of pattern restricted subsequences. Finally, similar results 
are also obtained for words generated by a Markov sequence.

\bigskip\bigskip

\section{Random permutations} \label{sec:rp}

The asymptotic behavior of the length of the longest alternating subsequence
has been studied by several authors, including Pemantle 
\cite[page~684]{stanley2006increasing}, 
Stanley~\cite{stanley511419longest}
and \newline Widom~\cite{widom2006limiting}, who by a mixture of generating function 
methods and saddle point techniques get the following result:

\begin{thm}
\label{th:altper}Let $\boldsymbol{\tau}$,  be a uniform random permutation in
the symmetric group $\mathcal{S}_{n}$,
and let ${\rm LA}_n(\tau)$ be the length of the longest alternating
subsequence of $\boldsymbol{\tau}$. Then, 
\begin{align*}
\mathbf{E}  \operatorname*{LA}\nolimits_{n}\left(  \boldsymbol{\tau
}\right)      & =\frac{2n}{3}+\frac{1}{6}\, ,\qquad n\ge 2\\
\operatorname*{Var} \operatorname*{LA}\nolimits_{n}\left(
\boldsymbol{\tau}\right)      & =\frac{8n}{45}-\frac{13}{180}\,,
\quad n\ge 4.
\end{align*}

Moreover, as $n\rightarrow\infty$,%
\[
\frac{\operatorname*{LA}_{n}(  \boldsymbol{\tau})  -{2n/3}%
\,}{\sqrt{\,8n/45}}\Longrightarrow\mathcal{Z}\text{,}%
\]
where $\mathcal{Z}$ is a standard normal random variable and where
$\Longrightarrow$ denotes convergence in distribution.
\end{thm}

The present section is devoted to give a simple probabilistic proof of the above
result. To provide 
 such a proof we make use of a well known correspondence which
transform the problem into that of counting the maxima of a sequence of iid
random variables uniformly distributed on $[0,1]$.
In order to establish the weak limit result, 
a central limit theorem for $m$-dependent random variables
is then briefly recalled.

\bigskip

Let us start by recalling some well known facts (Durrett \cite[Chapter~1]{durret},
Resnick \cite[Chapter~4]{resnick1999probability}). For each $n\geq1$
(including $n=\infty$), let $\mu_{n}$ be the uniform measure on $\left[
0,1\right]  ^{n}$ and, for each $n\geq1$, let the function $T_{n}%
:[0,1]^{n}\rightarrow\mathcal{S}_{n}$ be defined by $T_{n}(a_{1},a_{2}%
,\ldots,a_{n})=\tau^{-1}$, where $\tau$ is the unique permutation $\tau
\in\mathcal{S}_{n}$ that satisfies $a_{\tau_{1}}<a_{\tau_{2}}<\cdots
<a_{\tau_{n}}$. Note that $T_{n}$ is defined for all $a\in\left[  0,1\right]
^{n}$ except for those for which $a_{i}=a_{j}$ for some $i\neq j$, and this
set has $\mu_{n}$-measure zero. A well known fact, sometimes 
attributed to R\'{e}nyi
 \cite{resnick1999probability}, asserts that the pushforward
measure $T_{n}\mu_{n}$,
i.e., the image of $\mu_n$ by $T_n$,
 corresponds to the uniform measure on $\mathcal{S}_{n}%
$, which we denote by $\nu_{n}$. The importance of this fact relies in the
observation that the map $T_{n}$ is order preserving, that is, $a_{i}<a_{j}$
if and only if $\left(  T_{n}a\right)  _{i}<\left(  T_{n}a\right)  _{j}$. This
implies that any event in $\mathcal{S}_{n}$ has a canonical representative in
$\left[  0,1\right]  ^{n}$ in terms of the order relation of its components.
Explicitly, if we consider the language $L$ of the formulas with no
quantifiers, one variable, say $x$, and with atoms of the form $x_{i}<x_{j}%
$, $i,j\in\left[  n\right]  $, then any event of the form $\left\{
x:\varphi\left(  x\right)  \right\}  $ where $\varphi\in L$, has the same
probability in $\left[  0,1\right]  ^{n}$ and in $\mathcal{S}_{n}$ under the
uniform measure. To give some examples, events like $\left\{  x:x\text{ has an
increasing subsequence of length }k\right\}  $, $\left\{  x:x\text{ avoids the
permutation }\sigma\right\}  $, $\left\{  x:x\text{ has an alternating
subsequence of length }k\right\}  $ have the same probability in $\left[
0,1\right]  ^{n}$ and $\mathcal{S}_{n}$. In particular, it should be clear that
\begin{equation}
\operatorname*{LA}\nolimits_{n}\left(  \boldsymbol{\tau}\right)
\overset{d}{=}\operatorname*{LA}\nolimits_{n}%
(\boldsymbol{a})\text{,}\label{Eq:equaldist}%
\end{equation}
where $\boldsymbol{\tau}$ is a uniform random permutation in $\mathcal{S}_{n}$, $\boldsymbol{a}$ is a uniform random sequence in $\left[  0,1\right]
^{n}$ and where $d$ means equality in distribution.

\bigskip

\emph{Maxima and minima.} Next, we say that the sequence $a=(a_{1},a_{2}%
,\ldots,a_{n})$ has a local maximum at the index $k$ if (i)
$a_{k}>a_{k+1}$ or $k=n$, and (ii) $a_{k}>a_{k-1}$ or $k=1$. Similarly, we say
that $a$ has a local minimum at the index $k$ if (i) $a_{k}<a_{k+1}$
or $k=n$, and \linebreak(ii) $a_{k}<a_{k-1}$. An observation that comes in handy is the
fact that counting the length of the longest alternating subsequence is
equivalent to counting maxima and minima of the sequence (starting with a
local minimum). This is attributed to B\'{o}na in Stanley 
\cite{stanley2006increasing}; for completeness, we prove it next.

\begin{prop}
\label{thm:maxima}For $\mu_{n}$-almost all sequences $a=(a_{1},a_{2}%
,\ldots,a_{n})\in\lbrack0,1]^{n}$,
\begin{align}
\operatorname*{LA}\nolimits_{n}(a)  & =\#\text{ local maxima of }a\text{
}+\#\text{ local minima of }a\label{eq:maxima}\\
& =\boldsymbol{1}\left(  a_{n}>a_{n-1}\right)  +2\,\boldsymbol{1}\left(
a_{1}>a_{2}\right)  +2
{\textstyle\sum\limits_{k=2}^{n-1}}
\boldsymbol{1}\left(  a_{k-1}<a_{k}>a_{k+1}\right) .\label{eq:maxima2}%
\end{align}
\end{prop}

\begin{pf}
For $\mu_{n}$-almost all $a\in\lbrack0,1]^{n}$, $a_{i}\neq a_{j}$ whenever
$i\neq j$, therefore we can assume that $a$ has no repeated components. Let
$t_{1},\ldots,t_{r}$ be the positions, in increasing order, of the local
maxima of the sequence $a$, and let $s_{1},\ldots,s_{r^{\prime}}$ be the
positions, in increasing order, of the local minima of $a$, not including the
local minima before the position $t_{1}$. Notice that the maxima and minima
are alternating, that is, $t_{i}<s_{i}<t_{i+1}$ for every $i$, implying that
$r^{\prime}=r$ or $r^{\prime}=r-1$. Also notice, that in  case $r^{\prime
}=r-1$, necessarily $t_{r}=n$. Therefore, since $\left(  a_{t_{1}},a_{s_{1}%
},a_{t_{2}},a_{s_{2}},\ldots\right)  $ is an alternating subsequence of $a$,
we have $\operatorname*{LA}_{n}(a)\geq r+r^{\prime}=\#$ local maxima $+\#$
local minima.

To establish the opposite inequality, take a maximal sequence of indices
$\{\ell_{i}\}_{i=1}^{m}$ such that $\left(  a_{\ell_{i}}\right)  _{i=1}^{m}$
is alternating. Move every odd index upward, following the gradient of $a$
(the direction, left or right, in which the sequence $a$ increases), till it
reaches a local maximum of $a$. Next, move every even index downward,
following the gradient of $a$ (the direction, left or right, in which the
sequence $a$ decreases), till it reaches a local minimum of $a$. Notice,
importantly, that this sequence of motions preserves the order relation
between the indices, therefore the resulting sequence of indices $\{\ell
_{i}^{\prime}\}_{i=1}^{m}$ is still increasing and, in addition, it is a
subsequence of $\left(  t_{1},s_{1},t_{2},s_{2},\ldots\right)  $. Now, since
the sequence $\left(  a_{\ell_{i}^{\prime}}\right)  _{i=1}^{m}$ is
alternating, it follows that $LA_{n}(a)\leq\#$ local maxima $+\#$ local minima. Finally, associating every local maxima not in the $n-$th position
with the closest local minima to its right, we obtain a one to one
correspondence, which leads to 
(\ref{eq:maxima2}).
\end{pf}

\bigskip

\emph{Mean and variance. }The above correspondence 
allows us to easily compute the mean and the 
variance of the length of the longest
alternating subsequence by going `back and forth' between $\left[  0,1\right]
^{n}$ and $\mathcal{S}_{n}$. For instance, given a random uniform sequence
$\boldsymbol{a}=\left(  \boldsymbol{a}_{1},\ldots,\boldsymbol{a}_{n}\right)
\in\left[  0,1\right]  ^{n}$, let $M_{k}:=\boldsymbol{1}(  \boldsymbol{a}$
 has a local maximum at the index $k)$,
$k\in\left\{  2,\ldots,n-1\right\}  $. Then
\[
\mathbf{E}  M_{k}  =\mu_{n}(a_{k-1}<a_{k}>a_{k+1})=\mu_{3}%
(a_{1}<a_{2}>a_{3})=\nu_{3}(\tau_{1}<\tau_{2}>\tau_{3})\text{,}%
\]
where again, $\nu_{n}$ is the uniform measure on $\mathcal{S}_{n}$, $n\geq1$. 
The event, $\left\{  \tau_{1}<\tau_{2}>\tau_{3}\right\}  $ corresponds to the
permutations $\left\{  132,231\right\}  $, which shows that $\mathbf{E}
M_{k}  =1/3$. 

Similarly,
\[
\mathbf{E}  M_{1}  =\nu_{2}(\tau_{1}>\tau_{2})=1/2\text{ and
}\mathbf{E}  M_{n}  =\nu_{2}(\tau_{1}<\tau_{2})=1/2\text{.}%
\]
Plugging these values into (\ref{eq:maxima2}), we get that 
$$\mathbf{E}%
\operatorname*{LA}\nolimits_{n}(\boldsymbol{\tau})=\frac{2n}{3}+\frac{1}{6}.$$

To compute the variance of LA$_n(\tau)$, first note that $\operatorname*{Cov}\left(
M_{k},M_{k+r}\right)  =0$ whenever $r\geq3$, and that $\mathbf{E}\left[
M_{k}M_{k+1}\right]  =0$. Now, going again back and forth between $\left[
0,1\right]  ^{n}$ and $\mathcal{S}_{n}$, we also obtain

\begin{align*}
\mathbf{E}\left[  M_{k}M_{k+2}\right]    & =\nu_{5}(\tau_{1}<\tau_{2}>\tau
_{3}<\tau_{4}>\tau_{5})=2/15\text{,}\\
\mathbf{E}\left[  M_{1}M_{3}\right]    & =\nu_{4}(\tau_{1}>\tau_{2}<\tau
_{3}>\tau_{4})=1/6\text{ }%
\end{align*}

and 
\[
\mathbf{E}\left[  M_{n-2}M_{n}\right]  =\nu_{4}(\tau_{1}<\tau
_{2}>\tau_{3}<\tau_{4})=1/6\text{.}%
\]
This implies from Proposition~\ref{thm:maxima} and (\ref{Eq:equaldist}), that
$$\operatorname*{Var}\operatorname*{LA}\nolimits_{n}(\boldsymbol{\tau
})=\frac{8n}{45}-\frac{13}{180}.$$

\bigskip

\emph{Asymptotic normality.} Recall that collection of random variables $\left\{
X_{i}\right\}  _{i=1}^{\infty}$ is called $m$-dependent if $X_{t+m+1}$ is
independent of $\left\{  X_{i}\right\}  _{i=1}^{t}$ for every $t\geq1$. For
such sequences the strong law of large numbers extends in a straightforward
manner just partitioning the summand in appropriate sums of independent random
variables, but the extension of the central limit theorem to this context is
less trivial (although a `small block' - `big block' argument will do the
job). For this purpose recall also the following particular case of a
theorem due to Hoeffding and
Robbins \cite{hoeffding1948central} (which can be also found in standard texts
such as Durrett \cite[Chapter~7]{durret} or Resnick 
\cite[Chapter~8]{resnick1999probability}).

\begin{thm}
\label{REFERENCE}\ Let $\left(  X_{i}\right)_{i\geq1}$ be a sequence of
identical distributed $m$-dependent \emph{bounded} random variables. Then
\[
\frac{ X_{1}+\cdots +X_{n}  -n\mathbf{E}X_{1}}{\gamma\sqrt{n}%
}\Longrightarrow\mathcal{Z}\text{,}%
\]
where $\mathcal{Z}$ is a standard normal random variable, and where the variance
term is given by
\[
\gamma^{2}=\operatorname*{Var}X_{1}+2
{\textstyle\sum\limits_{t=2}^{m+1}}
\operatorname*{Cov}\left(  X_{1},X_{t}\right)  \text{.}%
\]

\end{thm}

Now, let $\boldsymbol{a}=\left(  \boldsymbol{a}_{1},\boldsymbol{a}_{2}%
,\ldots\right)  $ be a sequence of iid random variables uniformly distributed in $[0,1]$, and let $\boldsymbol{a}^{(n)}=(\boldsymbol{a}_{1},\ldots,\boldsymbol{a}_{n})$ be the restriction of
the sequence $\boldsymbol{a}$ to the first $n$ indices. Recalling
(\ref{Eq:equaldist}) and Proposition~\ref{thm:maxima}, it is clear that if
$\boldsymbol{\tau}$ is a uniform random permutation in $\mathcal{S}_{n}$,
\begin{equation}
\operatorname*{LA}\nolimits_{n}(\boldsymbol{\tau})\overset{d}%
{=}\boldsymbol{1}\left[  \boldsymbol{a}_{n}>\boldsymbol{a}_{n-1}\right]
+2\boldsymbol{1}\left[  \boldsymbol{a}_{1}>\boldsymbol{a}_{2}\right]  +2
{\textstyle\sum\limits_{k=2}^{n-1}}
\boldsymbol{1}\left[  \boldsymbol{a}_{k-1}<\boldsymbol{a}_{k}>\boldsymbol{a}
_{k+1}\right]  \text{,}\label{eq:fintoinf4}
\end{equation}
where $\overset{d}{=}$ denotes equality in distribution. Therefore,
since the random variables \linebreak$\left\{  \boldsymbol{1}\left[  \boldsymbol{a}
_{k-1}<\boldsymbol{a}_{k}>\boldsymbol{a}_{k+1}\right]  :k\geq2\right\}  $ are
identically distributed and $2$-dependent, we have by the strong law of large
numbers that with probability one 
\[
\lim_{n\rightarrow\infty}\frac{1}{n}
{\textstyle\sum\limits_{k=2}^{n-1}}
\boldsymbol{1}\left[  \boldsymbol{a}_{k-1}<\boldsymbol{a}_{k}>\boldsymbol{a}%
_{k+1}\right]  =\mu_{3}\left(  a_{1}<a_{2}>a_{3}\right)  =\frac{1}%
{3}\,.
\]
Therefore, from (\ref{eq:fintoinf4}) we get that, in probability,
\[
\lim_{n\rightarrow\infty}\frac{1}{n}\operatorname*{LA}\nolimits_{n}%
(\boldsymbol{\tau})=\frac{2}{3}\,.
\]

Finally, applying the above central limit theorem, we have
as $n\to\infty$
\begin{equation}
\frac{\operatorname*{LA}_{n}\left(  \boldsymbol{\tau}\right)  -2n/3}%
{\sqrt{\,n}\gamma}\Longrightarrow N(0,1),
\label{CLT}%
\end{equation}
where in our case, the variance term is given by%
\begin{align*}
\gamma^{2}  & =\operatorname*{Var}\left(  2\boldsymbol{1}\left[
\boldsymbol{a}_{1}<\boldsymbol{a}_{2}>\boldsymbol{a}_{3}\right]  \right)
+2\operatorname*{Cov}\left(  2\boldsymbol{1}\left[  \boldsymbol{a}%
_{1}<\boldsymbol{a}_{2}>\boldsymbol{a}_{3}\right]  ,2\boldsymbol{1}%
\left[  \boldsymbol{a}_{2}<\boldsymbol{a}_{3}>\boldsymbol{a}_{4}\right]
\right)  \\
&\quad  +2\operatorname*{Cov}\left(  2\boldsymbol{1}\left[  \boldsymbol{a}%
_{1}<\boldsymbol{a}_{2}>\boldsymbol{a}_{3}\right]  ,2\boldsymbol{1}%
\left[  \boldsymbol{a}_{3}<\boldsymbol{a}_{4}>\boldsymbol{a}_{5}\right]
\right)  \\
& =\frac{8}{45}\text{,}%
\end{align*}
from the computations carried out in the previous paragraph.

\bigskip

\begin{rem}
\rm The above approach via $m$-dependence has another advantage, it provides 
using standard m-dependent probabilistic statements various types of results 
on ${\rm LA}_n(\tau)$ such as, for example, the
exact fluctutation theory via the law of iterated logarithm.  In our
setting, it gives:
\begin{align*}
\limsup_{n\rightarrow\infty}\frac{\operatorname*{LA}\nolimits_{n}%
(\boldsymbol{\tau})-\mathbf{E}\,{\rm LA}_n(\tau)}{\sqrt{n\log\log n}}  
& =\frac{4}{3\sqrt{5}}\text{,}\\
\liminf_{n\rightarrow\infty}\frac{\operatorname*{LA}\nolimits_{n}%
(\boldsymbol{\tau})-\mathbf{E}\,{\rm LA}_n(\tau)}{\sqrt{n\log\log n}}  & =-\frac{4}{3\sqrt{5}}\text{.}%
\end{align*}
Besides the LIL, other types of probabilistic statements on LA$_n(\tau)$ are
possible, e.g., local limit theorems \cite{riauba1977local},
large deviations \cite{heinrich1985non}, exponential inequalities \cite{averkamp2005wavelet}, etc.  This types of statements are also true in the settings of our next sections.
\end{rem}

\bigskip

\section{Finite alphabet random words}

Consider a (finite) random sequence $\boldsymbol{a}=(\boldsymbol{a}_{1},\boldsymbol{a}%
_{2},\ldots,\boldsymbol{a}_{n})$ with distribution $\mu^{(n)}$, where $\mu$ is
a probability measure supported on a finite set $[q]=\{1,\ldots,q\}$. Our goal
now is to study the length of the longest alternating subsequence of the
random sequence $\boldsymbol{a}$. This new situation differs from the previous
one mainly in that the sequence can have repeated values. Thus, in order to
check if a point is a maximum or a minimum, it is not enough to `look at' its
nearest neighbors, losing the advantage of the $2$-dependence that we had in
the previous case.\ However, Instead, we can use the stationarity of the
property `being a local maximum' with respect to some extended sequence to
study the asymptotic behaviour of $\operatorname*{LA}_{n}\left(
\boldsymbol{a}\right)  $. As a matter of notation, we will use generically,
the expression $\operatorname*{LA}_{n}\left(  \mu\right)  $ for the
distribution of the length of the longest alternating subsequence of a
sequence $\boldsymbol{a}=(\boldsymbol{a}_{1},\boldsymbol{a}_{2},\ldots
,\boldsymbol{a}_{n})$ having the product distribution $\mu^{(n)}$.

In this section we proceed more or less along the lines of the previous
section, relating the counting of maxima to the length of the longest
alternating subsequence and then, through mixing and ergodicity, obtain
results on the asymptotic mean, variance, convergence of averages and
asymptotic normality of the longest alternating subsequence. These results are
presented in Theorem \ref{th:LLNiid} (convergence in probability), and Theorem
\ref{th:CLTiid} (asymptotic normality).

\bigskip

\emph{Counting maxima and minima.} Given a sequence $a=(a_{1},a_{2}%
,\ldots,a_{n})\in\left[  q\right]  ^{n}$, we say that $a$ has a local maximum
at the index $k$, if (i) $a_{k}>a_{k+1}$ or $k=n$, and if (ii) for some $j<k$,
$a_{j}<a_{j+1}=\cdots a_{k-1}=a_{k}$ or for all $j<k$, $a_{j}=a_{k}$.
Likewise, we say that $a$ has a local minimum at the index $k$, if (i)
$a_{k}<a_{k+1}$ or $k=n$, and if (ii) for some $j<k$, $a_{j}>a_{j+1}=\cdots
a_{k-1}=a_{k}$. The identity (\ref{eq:maxima}) can be generalized, in a
straightforward manner to this context, so that
\begin{align*}
\operatorname*{LA}\nolimits_{n}\left(  a\right)    & =\#\text{ local maxima of
}a+\#\text{ local minima of }a\\
& =\boldsymbol{1}\left(  a\text{ has a local maximum at }n\right)  +2
{\textstyle\sum\limits_{k=1}^{n-1}}
\boldsymbol{1}\left(  a\text{ has a local maximum at }k\right)  \text{.}%
\end{align*}
Now, the only difficulty in adapting the proof of Theorem \ref{thm:maxima} to
our current framework is when moving in the direction of the gradient when trying to
modify the alternating subsequence to consist of only maxima and minima. 
Indeed,
we could get stuck at an index of gradient zero that is
neither maximum nor minimum. But this difficulty can easily be overcome by
just deciding to move to the right whenever we get in such a situation. We then
end up with an alternating subsequence consisting of only maxima and minima
through order preserving moves.

\bigskip

\emph{Infinite bilateral sequences.} More generally, given an infinite
bilateral sequence \linebreak$a=\left(  \ldots,a_{-1},a_{0},a_{1},\ldots\right)
\in\left[  q\right]  ^{\mathbb{Z}}$, we say that $a$ has a local maximum at
the index $k$, if for some $j<k$, $a_{j}<a_{j+1}=\cdots=a_{k}>a_{k+1}$ and that $a$ has a local minimum at the index $k$, if for some $j<k$,
$a_{j}>a_{j+1}=\cdots=a_{k}<a_{k+1}$. Also, set $a^{\left(  n\right)
}=\left(  a_{1},\ldots,a_{n}\right)  $ to be the truncation of $a$ to the
first $n$ positive indices. An important observation is the following: Let
\[
A_{k}=\left\{  a\in\lbrack q]^{\mathbb{Z}}:\text{For some }j\leq0\text{,
}a_{j}>a_{j+1}=\cdots=a_{k}>a_{k+1}\right\},
\]
\[
A_{k}^{\prime}=\left\{  a\in\lbrack q]^{\mathbb{Z}}:\text{For some }%
j\leq0\text{, }a_{j}\neq a_{j+1}=\cdots=a_{k}\leq a_{k+1}\right\},
\]
and
\[
A_{k}^{\prime\prime}=\left\{  a\in\lbrack q]^{\mathbb{Z}}:\text{For some
}j\geq1\text{, }a_{j}<a_{j+1}=\cdots=a_{k}\leq a_{k+1}\right\}  \text{.}%
\]
Then, for any bilateral sequence $a\in\left[  q\right]  ^{\mathbb{Z}}$, we
have
\[
\boldsymbol{1}\left(  a^{\left(  n\right)  }\text{ has a local maximum at
}k\right)  =\boldsymbol{1}\left(  a\text{ has a local maximum at }k\right)
+\boldsymbol{1}_{A_{k}}\left(  a\right)  \text{, if }k<n\text{,}%
\]
and
\begin{align*}
\boldsymbol{1}\left(  a^{\left(  n\right)  }\text{ has a local maximum at
}n\right)    & =\boldsymbol{1}\left(  a\text{ has a local maximum at
}n\right)  \\
& +\boldsymbol{1}_{A_{n}}\left(  a\right)  +\boldsymbol{1}_{A_{n}^{\prime}%
}(a)+\boldsymbol{1}_{A_{n}^{\prime\prime}}(a)\text{.}%
\end{align*}
Hence,%
\begin{equation}
\operatorname*{LA}\nolimits_{n}(a^{(n)})=2
{\textstyle\sum_{k=1}^{n-1}}
\boldsymbol{1}\left(  a\text{ has a local maximum at }k\right)  +R_{n}\left(
a\right),  \label{eq:fintoinf3}%
\end{equation}
where the remainder term is given by 

\[
R_{n}\left(  a\right)  :=2{\sum\limits_{k=1}^{n-1}}\boldsymbol{1}_{A_{k}%
}(a)+\boldsymbol{1}\left(  a^{\left(  n\right)  }\text{ has a local maximum at
}n\right)  ,
\]
and is such that $\left\vert R_{n}\left(
a\right)  \right\vert \leq3$, since the sets $\left\{  A_{k}\right\}
_{k=1}^{n}$ are pairwise disjoint.

\bigskip

\emph{Stationarity.} Define the function $f:\left[  q\right]  ^{\mathbb{Z}%
}\rightarrow\mathbb{R}$ via 

\[
f\left(  a\right)  =2\,\boldsymbol{1}\left(  a\text{ has a local maximum at
the index }0\right)  .
\]

If $T:\left[  q\right]
^{\mathbb{Z}}\rightarrow\left[  q\right]  ^{\mathbb{Z}}$ is the (shift)
transformation such that $\left(  Ta\right)  _{i}=a_{i+1}$, and $T^{(k)}$ is
the $k$-th iterate of $T$, it is clear that $f\circ T^{(k)}(a)=2\,\boldsymbol{1}\left(  a\text{ has
a local maximum at }k\right)  $. With these notations,
(\ref{eq:fintoinf3}) becomes $\operatorname*{LA}\nolimits_{n}(a^{(n)})=
{\textstyle\sum\limits_{k=1}^{n-1}}
f\circ T^{(k)}(a)+R_{n}\left(  a\right)  $. In particular, if $\boldsymbol{a}$
is a random sequence with distribution $\mu^{\left(  \mathbb{Z}\right)  }$, 
and if
$T^{\left(  k\right)  }f$ is short for $f\circ T^{(k)}(\boldsymbol{a})$ the following
holds true:
\begin{equation}
\operatorname*{LA}\nolimits_{n}\left(  \mu\right)  \overset{d}{=}
{\textstyle\sum\limits_{k=1}^{n-1}}
T^{\left(  k\right)  }f+R_{n}\left(  \boldsymbol{a}\right)  \text{.}%
\label{eq:fintoinf2}%
\end{equation}

The transformation $T$ is measure preserving with respect to $\mu^{\left(
\mathbb{Z}\right)  }$ and, moreover, ergodic. Thus, by the classical
ergodic theorem (see, for example, \cite[Chapter~V]{shiryaev1996probability}), 
as $n\rightarrow\infty$, 
${\textstyle\sum\limits_{k=1}^{n}} T^{\left(  k\right)  }f/n\rightarrow\mathbf{E}f  $, where the
convergence occurs \emph{almost surely} and also \emph{in the mean}. The
limit can be easily computed:
\begin{align*}
\mathbf{E}f  & =2\,
{\textstyle\sum\limits_{k=0}^{\infty}}
\mathbb{P}\left(  \boldsymbol{a}_{-\left(  k+1\right)
}<\boldsymbol{a}_{-k}=\cdots=\boldsymbol{a}_{0}>\boldsymbol{a}_{1}\right) \\
& =2\,
{\textstyle\sum\limits_{k=0}^{\infty}}
{\textstyle\sum\limits_{x\in\left[  q\right]  }}
L_{x}^{2}p_{x}^{k+1}\\
& =2\,
{\textstyle\sum\limits_{x\in\left[  q\right]  }}
\frac{p_{x}}{1-p_{x}}L_{x}^{2}\\
& =
{\textstyle\sum\limits_{x\in\left[  q\right]  }}
\left(  \frac{L_{x}^{2}+U_{x}^{2}}{1-p_{x}}\right)  p_{x}\text{,}%
\end{align*}
where for $x\in\left[  q\right]  $, $p_{x}:=\mu\left(  \left\{
x\right\}  \right)  $, $L_{x}:=
{\textstyle\sum\limits_{y<x}}
p_{y}$ and $U_{x}:=
{\textstyle\sum\limits_{y>x}}
p_{y}$.

\bigskip

\emph{Oscillation.} Given a probability distribution $\mu$ supported on
$\left[  q\right]  $, define the `oscillation of $\mu$ at $x$', as $\operatorname*{osc}_{\mu}(x):={(L_{x}^{2}+U_{x}^{2})}/{(L_{x}+U_{x})
}$ and the total oscillation of the measure $\mu$ as $\operatorname*{Osc}%
\left(  \mu\right)  :=
{\textstyle\sum\limits_{x\in\lbrack q]}}
\operatorname*{osc}_{\mu}(x)p_{x}$. Interpreting the results of the previous
paragraph through  (\ref{eq:fintoinf2}), we conclude that

\begin{thm}
\label{th:LLNiid}Let $\boldsymbol{a}=(\boldsymbol{a}_{i})_{i=1}^{n}$ be a
sequence of iid random variables with common distribution $\mu$ supported
on $\left[  q\right]  $, and let $\operatorname*{LA}_{n}(\mu)$ be the length
of the longest alternating subsequence of $\boldsymbol{a}$. Then,
\[
\lim_{n\rightarrow\infty}\frac{\operatorname*{LA}\nolimits_{n}\left(
\mu\right)  }{n}=\operatorname*{Osc}\left(  \mu\right)  \text{, \emph{in the
mean}.}%
\]

\end{thm}

In particular, if $\mu$ a uniform distribution on $\left[  q\right]  $,
$\operatorname*{Osc}\left(  \mu\right)  =({2}/{3}-{1/3q})$, and thus
$\operatorname*{LA}_{n}\left(  \mu\right)/n $ is concentrated around
$({2/3}-{1/3q})$ both \emph{in the mean} and \emph{in probability}. We
should mention here that Mansour \cite{mansour2008longest}, using generating
function methods obtained, for $\mu$ uniform,
 an explicit formula for $\mathbf{E}%
\operatorname*{LA}_{n}\left(  \mu\right)  $, which, of
course, is asymptotically equivalent to $\left(  {2/3}-{1/3q}\right)  n$. From (\ref{eq:fintoinf2}) it is not difficult to
derive also a nonasymptotic expression for $\mathbf{E}\operatorname*{LA}_{n}\left(
\mu\right)  $:%
\begin{equation}
\mathbf{E}\operatorname*{LA}\nolimits_{n}\left(  \mu\right)
=n\operatorname*{Osc}\left(  \mu\right)  +
{\textstyle\sum_{x\in\lbrack q]}}
R_{1}(x)p_{x}+
{\textstyle\sum_{x\in\lbrack q]}}
R_{2}(x)p_{x}^{n}\text{,}\label{eq:unifiid}%
\end{equation}
where the terms $R_{1}(x)$ and $R_{2}(x)$ are given by:
\[
R_{1}(x)=\frac{L_{x}}{L_{x}+U_{x}}+\frac{2L_{x}U_{x}}{\left(  L_{x}%
+U_{x}\right)  ^{2}}-\operatorname*{osc}\nolimits_{\mu}(x)\text{\quad
and\quad}R_{2}(x)=\frac{U_{x}}{L_{x}+U_{x}}-\frac{2L_{x}U_{x}}{\left(
L_{x}+U_{x}\right)  ^{2}}\text{.}%
\]
Applying (\ref{eq:unifiid}) in the uniform case recovers computations as given in \cite{mansour2008longest}.

As far as the asymptotic limit of $\operatorname*{Osc}\left(  \mu\right)  $ is
concerned, we have the following bounds for a general $\mu$.

\begin{prop}
Let $\mu$ be a probability measure supported on the finite set $[q]$, then
\begin{equation}
\frac{1}{2}\left(  1-
{\textstyle\sum\limits_{x\in\lbrack q]}}
p_{x}^{2}\right)  \leq\operatorname*{Osc}\left(  \mu\right)  \leq\frac{2}%
{3}\left(  1-
{\textstyle\sum\limits_{x\in\lbrack q]}}
p_{x}^{3}\right)  \text{.}\label{bounds}%
\end{equation}

\end{prop}

\begin{pf}
Note that $
{\textstyle\sum\limits_{x\in\lbrack q]}}
L_{x}p_{x}=
{\textstyle\sum\limits_{i<j}}
p_{i}p_{j}=
{\textstyle\sum\limits_{x\in\lbrack q]}}
U_{x}p_{x}$ and $
{\textstyle\sum\limits_{x\in\lbrack q]}}
L_{x}p_{x}+
{\textstyle\sum\limits_{x\in\lbrack q]}}
U_{x}p_{x}+
{\textstyle\sum\limits_{x\in\lbrack q]}}
p_{x}^{2}=1$, which implies that
\begin{equation}
{\textstyle\sum\limits_{x\in\lbrack q]}}
L_{x}p_{x}=
{\textstyle\sum\limits_{x\in\lbrack q]}}
U_{x}p_{x}=\frac{1}{2}\left(  1-
{\textstyle\sum\limits_{x\in\lbrack q]}}
p_{x}^{2}\right)  \text{.}\label{eqvar1}%
\end{equation}
Similarly, for any permutation $\sigma\in S_{3}$, we have that $
{\textstyle\sum\limits_{x\in\lbrack q]}}
L_{x}U_{x}p_{x}=
{\textstyle\sum\limits_{i_{1}<i_{2}<i_{3}}}
p_{i_{1}}p_{i_{2}}p_{i_{3}}=
{\textstyle\sum\limits_{i_{\sigma\left(  1\right)  }<i_{\sigma\left(
2\right)  }<i_{\sigma\left(  3\right)  }}}
p_{i_{1}}p_{i_{2}}p_{i_{3}}$, which implies that $%
6 {\textstyle\sum\limits_{x\in\lbrack q]}}
L_{x}U_{x}p_{x}=
{\textstyle\sum\limits_{i_{1}\neq i_{2}\neq i_{3}}}
p_{i_{1}}p_{i_{2}}p_{i_{3}}$. Finally, an inclusion-exclusion argument leads
to%
\[
{\textstyle\sum\limits_{i_{1}\neq i_{2}\neq i_{3}}}
p_{i_{1}}p_{i_{2}}p_{i_{3}}=1-3
{\textstyle\sum\limits_{i_{i}=i_{2}}}
p_{i_{1}}p_{i_{2}}+2
{\textstyle\sum\limits_{i_{i}=i_{2}}}
p_{i_{1}}p_{i_{2}}p_{i_{3}}=1-3
{\textstyle\sum\limits_{x\in\left[  q\right]  }}
p_{x}^{2}+2
{\textstyle\sum\limits_{x\in\left[  q\right]  }}
p_{x}^{3}\,,
\]
and therefore
\begin{equation}
{\textstyle\sum\limits_{x\in\lbrack q]}}
L_{x}U_{x}p_{x}=\frac{1}{6}-\frac{1}{2}
{\textstyle\sum\limits_{x\in\left[  q\right]  }}
p_{x}^{2}+\frac{1}{3}
{\textstyle\sum\limits_{x\in\left[  q\right]  }}
p_{x}^{3}\text{.}\label{qvar2}%
\end{equation}

Now, to obtain the upper bound in (\ref{bounds}), note that
\begin{equation}
\operatorname*{Osc}\left(  \mu\right)  =
{\textstyle\sum\limits_{x\in\lbrack q]}}
\frac{L_{x}^{2}+U_{x}^{2}}{L_{x}+U_{x}}p_{x}=
{\textstyle\sum\limits_{x\in\lbrack q]}}
\left(  L_{x}+U_{x}\right)  p_{x}-2
{\textstyle\sum\limits_{x\in\lbrack q]}}
\frac{L_{x}U_{x}}{L_{x}+U_{x}}p_{x}\label{qvar3}%
\end{equation}
so that in particular, $\operatorname*{Osc}\left(  \mu\right)  \leq
{\textstyle\sum\limits_{x\in\lbrack q]}}
\left(  L_{x}+U_{x}\right)  p_{x}-2
{\textstyle\sum\limits_{x\in\lbrack q]}}
L_{x}U_{x}p_{x}$. Hence, using (\ref{eqvar1}) and (\ref{qvar2}),
$$\operatorname*{Osc}\left(  \mu\right)  \leq\frac{2}{3}\left(  1-
{\textstyle\sum\limits_{x\in\lbrack q]}}
p_{x}^{3}\right)  .$$

For the lower bound, note that $4
{\textstyle\sum\limits_{x\in\lbrack q]}}
\frac{L_{x}U_{x}}{L_{x}+U_{x}}p_{x}\leq
{\textstyle\sum\limits_{x\in\lbrack q]}}
\left(  L_{x}+U_{x}\right)  p_{x}$, and from (\ref{qvar3}) we get
\[
\operatorname*{Osc}\left(  \mu\right)  \geq\frac{1}{2}
{\textstyle\sum\limits_{x\in\lbrack q]}}
\left(  L_{x}+U_{x}\right)  p_{x}=\frac{1}{2}\left(  1-
{\textstyle\sum\limits_{x\in\lbrack q]}}
p_{x}^{2}\right)  \text{.}%
\]

\end{pf}

An interesting problem would be to determine the distribution $\mu$ over
$\left[  q\right]  $ that maximizes the oscillation. It is not hard to prove
that such an optimal distribution should be symmetric about $\left(
q-1\right)  /2$, but it is harder to establish its shape (at least
asymptotically in $q$). 

\bigskip

\emph{Mixing.} The use of ergodic properties to analyze the random variable
$\operatorname*{LA}_{n}\left(  \mu\right)  $ goes beyond the mere application
of the ergodic theorem. Indeed, the random variables $\left\{  T^{\left(
k\right)  }f:k\in\mathbb{Z}\right\}  $ introduced above exhibit mixing, or ``long range independence'', meaning that as $n\to\infty$
\[
\sup_{A\in\mathcal{F}_{\geq0},B\in\mathcal{F}_{<-n}}\left\vert
\mathbb{P}\left(  A\left\vert B\right.  \right)
-\mathbb{P}\left(  A\right)  \right\vert \rightarrow0,
\]
where, for $n\geq0$, $\mathcal{F}_{\geq n}$ (respectively $\mathcal{F}_{<n}$) is the $\sigma
$-field of events generated by \newline$\left\{  T^{\left(  k\right)  }f:k\geq
n\right\}  $ (respectively $\left\{  T^{\left(  k\right)  }f:k<n\right\}  $).
This kind of mixing condition is usually called \emph{uniformly strong mixing}
or $\varphi$-mixing , and the decreasing sequence
\begin{equation}
\varphi\left(  n\right)  :=\sup_{A\in\mathcal{F}_{\geq0},B\in\mathcal{F}%
_{<-n}}\left\vert \mathbb{P}\left(  A\left\vert B\right.  \right)
-\mathbb{P}\left(  A\right)  \right\vert, \label{eq:unifmixing}%
\end{equation}
is called the \emph{rate} of uniformly strong mixing (see, for
example, \cite[Chapter~1]{lin1996limit}). Below, Proposition
\ref{pro:vandependence} asserts that, in our case, such a rate decreases
exponentially. Let us prove the following lemma first.

\begin{lem}
\label{lemita}Let $\boldsymbol{a}=(\boldsymbol{a}_{i})_{i\in\mathbb{Z}}$ be a
bilateral sequence of iid random variables with common distribution $\mu$
supported on $\left[  q\right]  $. Let $C_{n,t}=\left\{  \boldsymbol{a}_{-n}=\cdots=\boldsymbol{a}_{-n+t-1}\neq
\boldsymbol{a}_{-n+t}\right\}  $, $n\geq 1$, \newline$0\leq t\leq n$, then:
\begin{enumerate}
\item[(i)] For any $A\in\mathcal{F}_{\geq0}$ and any $t\leq n$, the event
$C_{n,t}\cap A$ is independent of the $\sigma$-field $\mathcal{G}%
_{<-n}$ of events generated by $\left\{  \boldsymbol{a}_{i}:i<-n\right\}  $.

\item[(ii)] Restricted to the event $C_{n,t}$, the $\sigma$-fields
$\mathcal{F}_{\geq0}$ and $\mathcal{G}_{<-n}$ are independent.
\end{enumerate}
\end{lem}

\begin{pf}
Let the event $B_{r,s}:=\{\boldsymbol{a}_{r}<\boldsymbol{a}_{r+1}%
=\cdots=\boldsymbol{a}_{s}>\boldsymbol{a}_{s+1}\}$. Then, for $s_{1}%
<s_{2}<\cdots<s_{m}$, 
$\textstyle\prod_{i=1}^{m} T^(s_{i})f=\textstyle\sum {\textstyle\prod_{i=1}^{n}}
\boldsymbol{1}_{B_{r_{i},s_{i}}}$ holds true, where the sum runs over the $r_{1},\ldots
,r_{n}$ such that $s_{i-1}<r_{i}<s_{i}$ (letting $s_{0}=-\infty$) and where 
\[
f\left(  a\right)  =2\,\boldsymbol{1}\left(
a\text{ has a local maximum at the index }0\right)  .
\]
Now, since
the random variables $\left\{  T^{\left(  i\right)  }f,i\in\mathbb{Z}\right\}
$ are binary, then for any $A\in\mathcal{F}_{\geq0}$ the random
variable $\boldsymbol{1}_{A}$ can be expressed as a linear combination of terms of the
form $
{\textstyle\prod\limits_{i=1}^{m}}
T^{\left(  s_{i}\right)  }f$, where $0\leq s_{1}<\cdots<s_{m}$.

Next, $\boldsymbol{1}_{C_{n,t}}
{\textstyle\prod\limits_{i=1}^{m}}
T^{\left(  s_{i}\right)  }f=\boldsymbol{1}_{C_{n,t}}\left(
{\textstyle\sum}
{\textstyle\prod\limits_{i=1}^{n}}
\boldsymbol{1}_{B_{r_{i},s_{i}}}\right)  =\boldsymbol{1}_{C_{n,t}}\left(
{\textstyle\sum_{r_{1}\geq-n+t-1}}
{\textstyle\prod\limits_{i=1}^{n}}
\boldsymbol{1}_{B_{r_{i},s_{i}}}\right)  $, which implies that $\boldsymbol{1}_{C_{n,t}}%
{\textstyle\prod\limits_{i=1}^{m}}
T^{\left(  s_{i}\right)  }f$ and $\mathcal{G}_{<-n}$ are independent. This
implies, in particular, the independence of the events$\ C_{n,t}\cap A$ and
$B$, for any $A\in\mathcal{F}_{\geq0}$ and $B\in\mathcal{G}_{<-n}$, proving
(i). The statement (ii) follows directly from (i).
\end{pf}

\begin{prop}
\label{pro:vandependence}Let $\boldsymbol{a}=(\boldsymbol{a}_{i}%
)_{i\in\mathbb{Z}}$ be a bilateral sequence of iid random variables with
 $\mu$ supported on $\left[  q\right]  $. If the event $A$
belongs to the $\sigma$-field $\mathcal{F}_{\geq 0}$, then for any $n\geq 1$,
\[
\left\Vert \mathbb{P}\left(  A|\mathcal{G}_{<-n}\right)
-\mathbb{P}(A)\right\Vert _{\infty}
:=\sup\limits_{B\in
\mathcal{G}_{<-n}}\left\vert \mathbb{P}\left(  A\left\vert B\right.
\right)  -\mathbb{P}(A)\right\vert %
\leq2q\kappa^{n}\text{,}%
\]
where $\kappa:=\max\limits_{x\in\left[  q\right]  }\mu\left(  \left\{
x\right\}  \right)  $. In particular, the rate of uniform strong mixing of
the sequence $\left\{  T^{\left(  k\right)  }f:k\in\mathbb{Z}\right\}  $ (see
\eqref{eq:unifmixing}), satisfies $\varphi\left(  n\right)  \leq2q\kappa
^{n-1}$.
\end{prop}

\begin{pf}
Let $A\in\mathcal{F}_{\geq 0}$. By Lemma \ref{lemita}, 
$\mathbb{P}\left(  A\cap C_{n,r}\left\vert \mathcal{G}_{<-n}\right.
\right)  =\mathbb{P}\left(  A\cap C_{n,r}\right)  $, whenever $r\leq
n$. Therefore,
\begin{align*}
\mathbb{P}\left(  A\left\vert \mathcal{G}_{<-n}\right.  \right)   &
=
{\textstyle\sum\limits_{r=1}^{n}}
\mathbb{P}\left(  A\cap C_{n,r}\left\vert \mathcal{G}_{<-n}\right.
\right)  +\mathbb{P}\left(  A\cap\{\boldsymbol{a}_{-n}%
=\cdots=\boldsymbol{a}_{0}\}\left\vert \mathcal{G}_{<-n}\right.  \right) \\
&  =
{\textstyle\sum\limits_{r=1}^{n}}
\mathbb{P}\left(  A\cap C_{n,r}\right)  +\mathbb{P}\left(
A\cap\{\boldsymbol{a}_{-n}=\cdots=\boldsymbol{a}_{0}\}\left\vert
\mathcal{G}_{<-n}\right.  \right) \\
&  =\mathbb{P}(A)+\left(  \mathbb{P}\left(  A\cap
\{\boldsymbol{a}_{-n}=\cdots=\boldsymbol{a}_{0}\}\left\vert \mathcal{G}%
_{<-n}\right.  \right)  -\mathbb{P}\left(  A\cap\{\boldsymbol{a}%
_{-n}=\cdots=\boldsymbol{a}_{0}\}\right)  \right)  \text{.}%
\end{align*}
Then, it follows:%

\begin{align*}
\left\Vert \mathbb{P}\left(  A\left\vert \mathcal{G}_{<-n}\right.  \right)
-\mathbb{P}(A)\right\Vert _{\infty} &  \leq\mathbb{P}\left(  A\cap
\{\boldsymbol{a}_{-n}=\cdots=\boldsymbol{a}_{0}\}\right)  \\
&  +\left\Vert \mathbb{P}\left(  A\cap\{\boldsymbol{a}_{-n}=\cdots
=\boldsymbol{a}_{0}\}\left\vert \mathcal{G}_{<-n}\right.  \right)  \right\Vert
_{\infty}\\
&  \leq2\left\Vert \mathbb{P}\left(  \boldsymbol{a}_{-n}=\cdots=\boldsymbol{a}%
_{0}\left\vert \mathcal{G}_{<-n}\right.  \right)  \right\Vert _{\infty}\\
&  \leq2q\kappa^{n}%
\end{align*}

where the last conclusion follows trivially from $\mathcal{G}_{<-n}%
\supseteq\mathcal{F}_{\leq-\left(  n+1\right)  }$.
\end{pf}

Taking advantage of the mixing property we can now infer 
without much effort the behaviour of
the asymptotic variance and also deduce the asymptotic normality of
the statistic $\operatorname*{LA}_{n}\left(  \mu\right)  $ . This is done in the
next two paragraphs.

\bigskip

\emph{Variance. }The computation of the variance of the sequence $S_{n}=
{\textstyle\sum\limits_{k=1}^{n}}
T^{\left(  k\right)  }f$ is straightforward. Indeed%
\begin{equation}
\operatorname*{Var}  S_{n}  =n\left[  \operatorname*{Cov}\left(
f,f\right)  +2
{\textstyle\sum\limits_{k=1}^{n-1}}
\operatorname*{Cov}\left(  f,T^{\left(  k\right)  }f\right)  \right]  -2
{\textstyle\sum\limits_{k=1}^{n-1}}
k\operatorname*{Cov}\left(  f,T^{\left(  k\right)  }f\right)  \text{,}%
\label{covformula}%
\end{equation}
and the mixing property from Proposition \ref{pro:vandependence} implies that
$\left\vert \operatorname*{Cov}\left(  f,T^{\left(  k\right)  }f\right)
\right\vert $ decreases geometrically in $k$, so that all the series involved
in (\ref{covformula}) converge. Therefore,
\begin{equation}
\operatorname*{Var} S_{n} =n\gamma^{2}+\operatorname*{O}\left(
1\right)  \text{ ,\quad where }\quad\gamma^{2}=\operatorname*{Cov}\left(
f,f\right)  +2
{\textstyle\sum\limits_{k=1}^{n-1}}
\operatorname*{Cov}\left(  f,T^{\left(  k\right)  }f\right)  \text{.}%
\label{eq:variance1}%
\end{equation}

Moreover, for $k\leq l$, $\left\vert \operatorname*{Cov}\left(  \boldsymbol{1}%
_{A}(\boldsymbol{a}),T^{\left(  k\right)  }f\right)  \right\vert
\leq\mathbf{E}\boldsymbol{1}_{A}(\boldsymbol{a})\leq\kappa^{l}$, and for
$k\geq l$, and making use of Proposition \ref{pro:vandependence}, $\left\vert
\operatorname*{Cov}\left(  \boldsymbol{1}_{A}(\boldsymbol{a}),T^{\left(
k\right)  }f\right)  \right\vert \leq4q\kappa^{k-l-2}\mathbf{E}\boldsymbol{1}%
_{A}(\boldsymbol{a})\leq4q\kappa^{k-2}$. This implies that, as $n\rightarrow
\infty$,
\[
\left\vert \operatorname*{Cov}\left(
{\textstyle\sum\limits_{k=1}^{n-1}}
T^{\left(  k\right)  }f,
{\textstyle\sum\limits_{k=1}^{n-1}}
\boldsymbol{1}_{A_{k}}(\boldsymbol{a})\right)  \right\vert \leq4q^{3}\left\{
{\textstyle\sum\limits_{k\leq l}}
\kappa^{l}+
{\textstyle\sum\limits_{l\leq k}}
\kappa^{k}\right\}  =\operatorname*{O}\left(  1\right)  \text{.}%
\]
Similarly, using the Cauchy-Schwarz inequality, we have that
$\operatorname*{Cov}\left(
{\textstyle\sum\limits_{k=1}^{n}}
T^{\left(  k\right)  }f,\boldsymbol{1}_{\tilde{A}_{n}}(\boldsymbol{a})\right)
\rightarrow0$ where $\tilde{A}_{n}$ is either one of the events $A_{n}$,
$A_{n}^{\prime}$ or $A_{n}^{\prime\prime}$. Finally using the fact that
$\operatorname*{Cov}\left(
{\textstyle\sum\limits_{k=1}^{n-1}}
T^{\left(  k\right)  }f,T^{\left(  n\right)  }f\right)  =
{\textstyle\sum\limits_{k=1}^{n-1}}
\operatorname*{Cov}\left(  f,T^{\left(  k\right)  }f\right)  $ is bounded as
$n\rightarrow\infty$, we conclude that $\operatorname*{Cov}\left(
{\textstyle\sum\limits_{k=1}^{n-1}}
T^{\left(  k\right)  }f,R\left(  n\right)  \right)  =\operatorname*{O}\left(
1\right)  $, as $n\rightarrow\infty$. This implies the corresponding extension
of (\ref{eq:variance1}) to $\operatorname*{LA}_{n}\left(  \mu\right)  $:
\[
\operatorname*{Var}  \operatorname*{LA}\nolimits_	{n}(\mu)  =n\gamma^{2}+\operatorname*{O}\left(  1\right)  \text{ as
}n\rightarrow\infty\text{.}%
\]

Note that the bound just established is not
meaningless since the boundedness of $R_{n}\left(  \boldsymbol{a}\right)  $
only guarantees the weaker estimate $\operatorname*{Var}\operatorname*{LA}_n(\mu)   =n\gamma^{2}+\operatorname*{O}%
\left(  n^{1/2}\right)  $.

Let us proceed to compute $\gamma^{2}$. Let $f_{l}:\left[  q\right]  ^{\mathbb{Z}}\rightarrow\mathbb{R}$ via
\[
f_{l}\left(  a\right)  =2\,\boldsymbol{1}\left(
a_{-l}<a_{-l+1}=\cdots=a_{0}>a_{1}\right)  ,
\]
so that $f\left(  a\right)  =%
{\textstyle\sum\limits_{l=1}^{\infty}}
f_{l}\left(  a\right)  $. Note that
\[
\operatorname*{Cov}\left(  f,T^{\left(  k\right)  }f_{l}\right)  =\left\{
\begin{array}
[c]{cc}%
0 & \text{if }k\geq l+2\\
4%
{\textstyle\sum\limits_{x,y\in\left[  q\right]  }}
\left(  \frac{L_{x}}{1-p_{x}}\right)  \left(  L_{y}p_{y}^{l}\right)
L_{x\wedge y}p_{x}-2\operatorname*{Osc}\left(  \mu\right)
{\textstyle\sum\limits_{y\in\left[  q\right]  }}
L_{y}^{2}p_{y}^{l} & \text{if }k=l+1\\
-2\operatorname*{Osc}\left(  \mu\right)
{\textstyle\sum\limits_{y\in\left[  q\right]  }}
L_{y}^{2}p_{y}^{l} & \text{if }1\leq k\leq l\\
4%
{\textstyle\sum\limits_{y\in\left[  q\right]  }}
L_{y}^{2}p_{y}^{l}-2\operatorname*{Osc}\left(  \mu\right)
{\textstyle\sum\limits_{y\in\left[  q\right]  }}
L_{y}^{2}p_{y}^{l} & \text{if }0=k\leq l,
\end{array}
\right.  
\]
and thus
\begin{align*}
\gamma^{2}  & =\operatorname*{Var}f +2
{\textstyle\sum\limits_{k=1}^{\infty}}
{\textstyle\sum\limits_{l=k-1}^{\infty}}
\operatorname*{Cov}\left(  f,T^{\left(  k\right)  }f_{l}\right)  \\
& =\operatorname*{Osc}\left(  \mu\right)  \left(  2-3\operatorname*{Osc}
\left(  \mu\right)  -4
{\textstyle\sum\limits_{x\in\left[  q\right]  }}
\left(  \frac{L_{x}}{1-p_{x}}\right)  ^{2}p_{x}\right)  +8
{\textstyle\sum\limits_{x,y\in\left[  q\right]  }}
\frac{L_{x}L_{y}L_{x\wedge y}}{\left(  1-p_{x}\right)  \left(  1-p_{y}\right)
}p_{x}p_{y}\text{.}
\end{align*}

We further mention at this point that 
Mansour \cite{mansour2008longest}
already obtained, with generating function methods, an exact expression for the variance when $\mu$ is the
uniform distribution on $\left[  q\right]  $. It is given (as it can also be
checked from (\ref{eq:variance1})) by 
$$\gamma^{2}=\frac{8}{45}\left[
\frac{\left(  1+1/q\right)  (1-3/4q)(1-1/2q)}{\left(  1-1/2q\right)  }\right]
.$$

\bigskip{}

\emph{Asymptotic normality. }Under appropriate conditions (say, asymptotic
positive variance and fast enough mixing), it is natural to expect for the
sequence of partial sums to be asymptotically normal. In our model, 
this is indeed 
the case. Let us recall the following central limit theorem which goes back to
Volkonskii and Rozanov \cite[Theorem 1.2]{volkonskii1959some} and which can be
found, greatly generalized, in texts such as  Bradley \cite[Theorem 10.3]{bradley2007introduction}.

\begin{thm} \label{thm:cltmixing}
Let $\boldsymbol{x}=\left(  \boldsymbol{x}%
_{i}\right)  _{i\in\mathbb{Z}}$ be a strictly stationary sequence of bounded
random variables such that the sequence 
\[
\alpha(n){:=}
\sup\limits_{A\in\mathcal{F}_{\geq0},B\in\mathcal{F}_{<-n}}\left\vert
\mathbb{P}\left(  A\cap B\right)  -\mathbb{P}
(A)\mathbb{P}(B)\right\vert 
\]
is summable (i.e. $
{\textstyle\sum\limits_{n\geq1}}
\alpha\left(  n\right)  <\infty$), where $\mathcal{F}_{\geq0}$ is the $\sigma
$-field generated by the random variables $\left\{  \boldsymbol{x}%
_{i}:i\geq0\right\}  $ and $\mathcal{F}_{<-n}$, $n\geq1$, is the $\sigma$-field 
generated by the random variables $\left\{  \boldsymbol{x}%
_{i}:i<-n\right\}  $. Then,

\begin{enumerate}
\item $\gamma^{2}:=\operatorname*{Var} \boldsymbol{x}_{0}
+2
{\textstyle\sum\limits_{t=1}^{\infty}}
\operatorname*{Cov}(\boldsymbol{x}_{0},\boldsymbol{x}_{t})$ exists in
$[0,\infty)$, the sum being absolutely convergent.

\item If $\gamma^{2}>0$, then as $n\rightarrow\infty$,
$$\frac{
{\textstyle\sum\limits_{t=1}^{n}}
\boldsymbol{x}_{t}-n\mathbf{E}\boldsymbol{x}_{0}  }{\sqrt{n}\gamma
}\Longrightarrow\mathcal{Z},$$ 
where $\mathcal{Z}$ is a standard normal
random variable.
\end{enumerate}
\end{thm}

Now, the asymptotic normality of $\operatorname*{LA}_{n}\left(  \mu\right)  $,
namely, the fact that as $n\rightarrow\infty$, 
$$\frac{
\operatorname*{LA}_{n}(  \mu)  -n\operatorname*{Osc}(
\mu)    }{\sqrt{n}\gamma}\Longrightarrow\mathcal{Z},$$ 
is clear:
By Proposition \ref{pro:vandependence}, the mixing coefficients $\alpha(n)$
decrease geometrically, implying the summability of $
{\textstyle\sum}
\alpha\left(  n\right)  $. Summarizing, we get:

\begin{thm}
\label{th:CLTiid}Let $\boldsymbol{a}=(\boldsymbol{a}_{i})_{i=1}^{n}$ be a
sequence of iid random variables, with common distribution $\mu$ supported
on $\left[  q\right]  $, and let $\operatorname*{LA}_{n}(\mu)$ be the length
of the longest alternating subsequence of $\boldsymbol{a}$. Then, as
$n\rightarrow\infty$,
\[
\frac{ \operatorname*{LA}_{n}(  \mu)  -n\operatorname*{Osc}
(  \mu)  }{\sqrt{n}\gamma}\Longrightarrow\mathcal{Z}
\text{, }
\]
where $\mathcal{Z}$ is a standard normal random variable and $\gamma$ is
given by \eqref{eq:variance1}.
\end{thm}

\begin{rem} It is clear that the above proofs extend to countable infinite alphabets, without major modification. A parallel situation for the longest increasing subsequence is given in  \cite{houdré2008longest}, though in that context a more delicate ``sandwich'' argument is required.
\end{rem}
\bigskip

\section{Markovian words}
\label{sec:Markovian}

Our probabilistic methodologies also provide results beyond the iid framework. Let now $\left(  \boldsymbol{x}%
_{k}\right)  _{k\geq0}$ be an ergodic Markov chain started at stationarity and whose state space is a finite
linearly ordered set $\mathcal{A}$, so that without loss of generality,
$\mathcal{A}=\left[  q\right]  $. Our objective (as before), is to study the
behavior of the statistics $\operatorname*{LA}_{n}\left(  \boldsymbol{x}%
_{0},\ldots,\boldsymbol{x}_{n}\right)  $.\bigskip

\emph{Adding gradient information to the chain. }Let us consider the related
process $\left(  \boldsymbol{y}_{k}\right)  _{k\geq0}$ defined recursively as follows:

\begin{description}
\item[-] $\boldsymbol{y}$$_{0}=1$.

\item[-] $\boldsymbol{y}$$_{k+1}=1$ if $\boldsymbol{x}$$_{k+1}>$%
$\boldsymbol{x}$$_{k}$ or if $\boldsymbol{x}$$_{k+1}=$$\boldsymbol{x}$$_{k}$
and $\boldsymbol{y}$$_{k}=1$.

\item[-] $\boldsymbol{y}_{k+1}=-1$ if $\boldsymbol{x}_{k+1}<$
$\boldsymbol{x}_{k}$ or if $\boldsymbol{x}_{k+1}=\boldsymbol{x}_{k}$
and $\boldsymbol{y}_{k}=-1$.
\end{description}

This new sequence basically carries the information indicating that the sequence is
increasing or decreasing at $k$ (we define the sequence $\boldsymbol{x}%
_{1},\boldsymbol{x}_{2},\ldots$ to be increasing at $k$ if $\boldsymbol{x}%
_{k}>\boldsymbol{x}_{k-1}$ or if it is increasing at $k-1$ and $\boldsymbol{x}%
_{k}=\boldsymbol{x}_{k-1}$, likewise, the sequence is decreasing at $k$ if
$\boldsymbol{x}_{k}<\boldsymbol{x}_{k-1}$ or if it is decreasing at $k-1$ and
$\boldsymbol{x}_{k}=\boldsymbol{x}_{k-1}$).

The following holds true for the process $\left(  \boldsymbol{x}%
_{k},\boldsymbol{y}_{k}\right)  _{k\geq0}$:

\begin{prop}
\label{pro:stat}The process $\left(  \boldsymbol{x}_{k},\boldsymbol{y}%
_{k}\right)  _{k\geq0}$ is Markov, with transition probabilities given by
\begin{align*}
p_{\left(  r,\pm1\right)  \rightarrow\left(  s,1\right)  } &  =p_{r,s}
\boldsymbol{1}\left(  s>r\right)  \text{,\quad}p_{\left(  r,1\right)
\rightarrow\left(  r,1\right)  }=p_{r,r}\\
p_{\left(  r,\pm1\right)  \rightarrow\left(  s,-1\right)  } &  =p_{r,s}
\boldsymbol{1}\left(  s<r\right)  \text{,\quad}p_{\left(  r,-1\right)
\rightarrow\left(  r,-1\right)  }=p_{r,r}
\end{align*}
and stationary measure given by
\[
\pi_{\left(  r,1\right)  }=\left(  1-p_{r,r}\right)  ^{-1}
{\textstyle\sum\limits_{s<r}}
\pi_{s}p_{s,r}\text{,\quad}\pi_{\left(  r,-1\right)  }=\left(  1-p_{r,r}
\right)  ^{-1}
{\textstyle\sum\limits_{s>r}}
\pi_{s}p_{s,r}\text{.}
\]

Moreover, the Markov process $\left(  \boldsymbol{x}_{k},\boldsymbol{y}%
_{k-1},\boldsymbol{y}_{k}\right)  _{k\geq0}$ has a stationary measure given
by
\begin{align*}
\pi_{\left(  r,1,1\right)  } &  =
{\textstyle\sum\limits_{t<s\leq r}}
\frac{\pi_{t}p_{t,s}p_{s,r}}{1-p_{s,s}}\text{,\quad}\pi_{\left(
r,-1,-1\right)  }=
{\textstyle\sum\limits_{t>s\geq r}}
\frac{\pi_{t}p_{t,s}p_{s,r}}{1-p_{s,s}}\\
\pi_{\left(  r,1,-1\right)  } &  =
{\textstyle\sum\limits_{t<s>r}}
\frac{\pi_{t}p_{t,s}p_{s,r}}{1-p_{s,s}}\text{,\quad}\pi_{\left(
r,-1,1\right)  }=
{\textstyle\sum\limits_{t>s<r}}
\frac{\pi_{t}p_{t,s}p_{s,r}}{1-p_{s,s}}%
\end{align*}

\end{prop}

\begin{pf}
The process is Markov since by definition $\boldsymbol{y}_{k+1}\in
\sigma\left(  \boldsymbol{x}_{k},\boldsymbol{x}_{k+1},\boldsymbol{y}%
_{k}\right)  $ and since $\left(  \boldsymbol{x}_{k}\right)  _{k\geq0}$ is
Markov. The transition probabilities are easily obtained from the definition,
and moreover,
\begin{align*}
&
{\textstyle\sum\limits_{r}}
\pi_{\left(  r,1\right)  }p_{\left(  r,1\right)  \rightarrow\left(
u,1\right)  }+
{\textstyle\sum\limits_{r}}
\pi_{\left(  r,-1\right)  }p_{\left(  r,-1\right)  \rightarrow\left(
u,1\right)  }\\
&  =
{\textstyle\sum\limits_{r\leq u}}
\left(  1-p_{r,r}\right)  ^{-1}
{\textstyle\sum\limits_{t<r}}
\pi_{t}p_{t,r}p_{r,u}+
{\textstyle\sum\limits_{r<u}}
\left(  1-p_{r,r}\right)  ^{-1}
{\textstyle\sum\limits_{t>r}}
\pi_{t}p_{t,r}p_{r,u}\\
&  =
{\textstyle\sum\limits_{r<u}}
\left(  1-p_{r,r}\right)  ^{-1}
{\textstyle\sum\limits_{t\neq r}}
\pi_{t}p_{t,r}p_{r,u}+\left(  1-p_{u,u}\right)  ^{-1}
{\textstyle\sum\limits_{t<u}}
\pi_{t}p_{t,u}p_{u,u}\\
&  =
{\textstyle\sum\limits_{t<u}}
\pi_{t}p_{t,u}+\left(  1-p_{u,u}\right)  ^{-1}
{\textstyle\sum\limits_{t<u}}
\pi_{t}p_{t,u}p_{u,u}\\
&  =\pi_{\left(  u,1\right)  }\text{.}%
\end{align*}
Similar computations show that
\[
{\textstyle\sum\limits_{r}}
\pi_{\left(  r,1\right)  }p_{\left(  r,1\right)  \rightarrow\left(
u,-1\right)  }+
{\textstyle\sum\limits_{r}}
\pi_{\left(  r,-1\right)  }p_{\left(  r,-1\right)  \rightarrow\left(
u,-1\right)  }=\pi_{\left(  u,-1\right)  }\text{,}
\]
thus proving that $\pi_{\left(  u,\pm1\right)  }$ is the stationary measure of
$\left(  \boldsymbol{x}_{k},\boldsymbol{y}_{k}\right)  _{k\geq0}$.

For the chain $\left(  \boldsymbol{x}_{k},\boldsymbol{y}_{k-1},\boldsymbol{y}
_{k}\right)  _{k\geq1}$ let us only verify one case since the others are
similar:%
\begin{align*}
&
{\textstyle\sum\limits_{r}}
\pi_{\left(  r,1,1\right)  }p_{_{\left(  r,1,1\right)  }\rightarrow_{\left(
u,1,1\right)  }}+
{\textstyle\sum\limits_{r}}
\pi_{\left(  r,-1,1\right)  }p_{_{\left(  r,-1,1\right)  }\rightarrow_{\left(
u,1,1\right)  }}\\
&  =
{\textstyle\sum\limits_{r\leq u}}
{\textstyle\sum\limits_{t<s\leq r}}
\frac{\pi_{t}p_{t,s}p_{s,r}}{1-p_{s,s}}p_{r,u}+
{\textstyle\sum\limits_{r\leq u}}
{\textstyle\sum\limits_{t>s<r}}
\frac{\pi_{t}p_{t,s}p_{s,r}}{1-p_{s,s}}p_{r,u}\\
&  =
{\textstyle\sum\limits_{s<r\leq u}}
\frac{p_{s,r}}{1-p_{s,s}}p_{r,u}
{\textstyle\sum\limits_{t<s}}
\pi_{t}p_{t,s}+
{\textstyle\sum\limits_{s<r\leq u}}
\frac{p_{s,r}}{1-p_{s,s}}p_{r,u}
{\textstyle\sum\limits_{t>s}}
\pi_{t}p_{t,s}+
{\textstyle\sum\limits_{s=r\leq u}}
\frac{p_{s,r}}{1-p_{s,s}}p_{r,u}
{\textstyle\sum\limits_{t<s}}
\pi_{t}p_{t,s}\\
&  =
{\textstyle\sum\limits_{s<r\leq u}}
\pi_{s}p_{s,r}p_{r,u}+
{\textstyle\sum\limits_{s<r\leq u}}
\frac{\pi_{s}p_{s,r}p_{r,u}p_{r,r}}{1-p_{r,r}}\\
&  =\pi_{\left(  u,1,1\right)  }\text{.}%
\end{align*}

\end{pf}

\bigskip

\emph{Oscillations of a Markov chain. }Given an ergodic Markov chain
$\boldsymbol{x}:=\left(  \boldsymbol{x}_{k}\right)  _{k\geq1}$ whose state
space is a finite linearly ordered set, define 
\begin{align*}
\operatorname*{Osc}\nolimits^{+}\left(  \boldsymbol{x}\right)    &
:={\sum\limits_{t<s>r}}{(\pi_{t}p_{t,s}p_{s,r})}/{(1-p_{s,s})}\\
\operatorname*{Osc}\nolimits^{-}\left(  \boldsymbol{x}\right)    &
:={\sum\limits_{t>s<r}}(\pi_{t}p_{t,s}p_{s,r})/(1-p_{s,s})
\end{align*}
and $\operatorname*{Osc}\left(
\boldsymbol{x}\right)  :=\operatorname*{Osc}^{+}\left(  \boldsymbol{x}\right)
+\operatorname*{Osc}^{-}\left(  \boldsymbol{x}\right)  $ (
$=2\operatorname*{Osc}^{+}\left(  \boldsymbol{x}\right)  =2\operatorname*{Osc}%
^{-}\left(  \boldsymbol{x}\right)  $ ). With these notations, we have: 

\begin{thm}
Let $\operatorname*{LA}\nolimits_{n}\left(  \boldsymbol{x}_{0},\ldots
,\boldsymbol{x}_{n}\right)  $ be the length of the longest alternating
subsequence of the first $n+1$ elements of the Markov chain $\left(
\boldsymbol{x}_{k}\right)  _{k\geq0}$. Then, as $n\rightarrow\infty$,
\[
\frac{\operatorname*{LA}\nolimits_{n}\left(  \boldsymbol{x}_{0},\ldots
,\boldsymbol{x}_{n}\right)  }{n}\rightarrow\operatorname*{Osc}\left(
\boldsymbol{x}\right)  \text{,}%
\]
in the mean and almost surely.
\end{thm}

\begin{pf}
From the very definition of $\boldsymbol{y}_{k}$,
\[
\operatorname*{LA}\nolimits_{n}\left(  \boldsymbol{x}_{0},\ldots
,\boldsymbol{x}_{n}\right)  =
{\textstyle\sum\limits_{k=0}^{n-1}}
\boldsymbol{1}\left(  \boldsymbol{y}_{k}\boldsymbol{y}_{k+1}=-1\right)  \text{,}%
\]
therefore, by the ergodic theorem,%
\[
\frac{\operatorname*{LA}\nolimits_{n}\left(  \boldsymbol{x}_{0},\ldots
,\boldsymbol{x}_{n}\right)  }{n}\rightarrow \Pi \left(  \boldsymbol{y}_{0}\boldsymbol{y}_{1}=-1\right)  \text{,}%
\]
in the mean and almost surely and where $\Pi$ is the stationary measure of the chain. Now, from Proposition \ref{pro:stat},
\[
\Pi\left(  \boldsymbol{y}_{0}\boldsymbol{y}%
_{1}=-1\right)  =
{\textstyle\sum\limits_{t<s>r}}
\frac{\pi_{t}p_{t,s}p_{s,r}}{1-p_{s,s}}+
{\textstyle\sum\limits_{t>s<r}}
\frac{\pi_{t}p_{t,s}p_{s,r}}{1-p_{s,s}}\text{,}%
\]
from which the result follows.
\end{pf}

\begin{rem}
Above, the case $p_{t,s}=p_{s}$ (and therefore $\pi_{t}=p_{t}$), corresponds to iid
letters thus recovering Theorem \ref{th:LLNiid}.
\end{rem}

\bigskip

\emph{Central limit theorem: } In case the asymptotic variance term of
$\operatorname*{LA}\nolimits_{n}\left(  \boldsymbol{x}_{1},\ldots
,\boldsymbol{x}_{n}\right)  $ is nonzero, then since \\*
$\operatorname*{LA}\nolimits_{n}\left(  \boldsymbol{x}_{1},\ldots
,\boldsymbol{x}_{n}\right)  $ is an additive functional of the finite Markov
chain $\left(  \boldsymbol{x}_{k},\boldsymbol{y}_{k-1},\boldsymbol{y}%
_{k}\right)  _{k\geq0}$, and since the mixing rate of an ergodic Markov chain with finite state space is exponentially decreasing, Theorem \ref{thm:cltmixing} imply that, for some $\gamma>0$, 
\[
\frac{ \operatorname*{LA}_{n}\left(  \boldsymbol{x}_{0},\ldots
,\boldsymbol{x}_{n}\right)  -n\operatorname*{Osc}\left(  \boldsymbol{x}%
\right)    }{\sqrt{n}\gamma}\Longrightarrow\mathcal{Z}\text{, }%
\]
where $\mathcal{Z}$ is a standard normal random variable. The reader should
contrast this last fact with the increasing subsequence results where the iid and Markov limiting laws differ when the alphabet has a size of four or more (\cite{trevis2}).

\section{Concluding remarks}

Determining the length of the longest alternating subsequence of a random 
\emph{pattern-avoiding} permutation or word, has been recently studied by Firro, Mansour and Wilson \cite{firro2006three,firro2007longest,mansour2008blongest} inspired by the work of Deutsch, Hildebrand and Wilf \cite{deutsch51longest} on the longest increasing subsequence of pattern-avoiding permutations. In such a case, a probabilistic (i.e. measure theoretic) approach is also possible once an appropriate recursive description of the pattern-avoiding permutations is given. Such recursive description is the subject of an extensive list of works, originating from an old standing conjecture of Zeilberger \cite{zeilberger1990holonomic} claiming in particular, that the set of pattern avoiding permutations is $P$-recursive. In the case of avoiding patterns of length $3$ a concise work is found in \cite{firro2006three}. A canonical example of this situation is the case of permutations avoiding the pattern $(123)$, or equivalently, sequences in $[0,1]^n$ avoiding the pattern $(123)$ (recall the observation at the beginning of Section \ref{sec:rp}). In this context, if we let $\mathcal{G}_{n}$ to be the set of sequences in $\left[ 0,1\right] ^{n}$ that avoid the pattern $\left( 123\right) $, and for $n\geq 1$ let%
\[
\upsilon _{n}\left( x_{n},\ldots ,x_{1}\right) =dx_{n}\ldots dx_{1}%
\boldsymbol{1}\left( \left( x_{n},\ldots ,x_{1}\right) \in \mathcal{G}%
_{n}\right) \text{,}
\]%
then, the recursive construction $\upsilon _{1}=dx_{1}$ and

\begin{align*}
\upsilon_{n+1}\left(  x_{n+1},\ldots,x_{1}\right)    & =dx_{n+1}\upsilon
_{n}\left(  x_{n},\ldots,x_{1}\right)  \boldsymbol{1}\left(  x_{n+1}%
>x_{n}\right)  \\
& +dx_{n}\upsilon_{n}\left(  x_{n+1},x_{n-1},\ldots,x_{1}\right)
\boldsymbol{1}\left(  x_{n}>\max\left\{  x_{1},\ldots,x_{n-1},x_{n}\right\}
\right)  \text{.}%
\end{align*}

for $n\geq 1$, holds. This recursive formulation for the restricted measure translates to a recursive formula for the distribution of the number of local maxima of the sequence $\left( x_{n},\ldots ,x_{1}\right) $
on $\mathcal{G}_{n}$: Let $M_{n}=\max \left\{ x_{1},\ldots ,x_{n}\right\} $,
let $L_{n}=\#\{i:x_{i}<x_{i+1}>x_{i+2}$, $i=1,\ldots ,n-2\}$ and let $\chi
_{n}=\boldsymbol{1}\left( M_{n}=x_{n}\right) $, $\varrho _{n}=\boldsymbol{1}\left(
x_{n}<x_{n-1}>x_{n-2}\right) $, then%

\begin{align*}
& \upsilon_{n+1}\left(  M_{n+1}=m,x_{n+1}=x,L_{n}=k,\chi_{n+1}=0,\varrho
_{n+1}=1\right)  \\
& =\upsilon_{n}\left(  M_{n}<m,x_{n}=x,L_{n}=k,\chi_{n}=0,\varrho
_{n}=1\right)  dm\\
& +\upsilon_{n}\left(  M_{n}<m,x_{n}=x,L_{n}=k-1,\chi_{n}=0,\varrho
_{n}=0\right)  dm\\
& +\upsilon_{n}\left(  M_{n}=x,x_{n}=x,L_{n}=k-1,\chi_{n}=1,\varrho
_{n}=0\right)  dm\\
& \upsilon_{n+1}\left(  M_{n+1}=m,x_{n+1}=x,L_{n}=k,\chi_{n+1}=0,\varrho
_{n+1}=0\right)  \\
& =\upsilon_{n}\left(  M_{n}=m,x_{n}<x,L_{n}=k,\chi_{n}=0\right)  dx\\
& \upsilon_{n+1}\left(  M_{n+1}=x,x_{n+1}=x,L_{n}=k,\chi_{n+1}=1,\varrho
_{n+1}=0\right)  \\
& =\upsilon_{n}\left(  M_{n}<x,x_{n}<x,L_{n}=k,\chi_{n}=0\right)  dx\\
& +\upsilon_{n}\left(  M_{n}<x,x_{n}<x,L_{n}=k,\chi_{n}=1,\varrho
_{n}=0\right)  dx\text{.}%
\end{align*}

These formulas can be interpreted as Markovian formulations of the process
of counting local maxima (therefore, the length of the longest alternating
subsequence), in the restricted space of permutations avoiding the pattern $%
\left( 123\right)$. Therefore the appropriate extension of the methods of Section \ref{sec:Markovian} lead to the corresponding results in this context. Notice however, that such Markovian formulation is not measure preserving, and the corresponding modification of the ergodic theorem, central limit theorem, etc., should be introduced. It is our goal in subsequent research, to study these methods for tractable (in the above sense), sets of pattern avoiding permutations or words, following this alternative probabilistic path just presented.


\begin{thebibliography}{10}

\bibitem{averkamp2005wavelet}
{\sc Averkamp, R.} and {\sc Houdr{\'e}, C.}, ``{Wavelet thresholding for nonnecessarily Gaussian noise: functionality},'' {\em Annals of statistics}, vol.~33, no.~5,
  pp.~2164--2193, 2005.
  
\bibitem{bradley2007introduction}
{\sc Bradley, R.}, {\em {Introduction to strong mixing conditions}}.
\newblock Kendrick Press, Heber City, Utah, 2007.

\bibitem{deutsch51longest}
{\sc Deutsch, E.}, {\sc Hildebrand, A.J.} and {\sc Wilf, H.S.}, ``{Longest increasing subsequences in pattern-restricted permutations},'' {\em The electronic journal of combinatorics}, vol.~9(2), no.~R12, 2003.


\bibitem{durret}
{\sc Durrett, R.}, {\em {Probability: Theory and Examples}}.
\newblock Thomson, 2005.


\bibitem{firro2006three} \textsc{Firro, G.}, \textsc{Mansour, T.} and 
\textsc{Wilson, M.C.}, ``{Three-letter-pattern-avoiding permutations and functional equations},'' \emph{Elect. J. Combin.},
vol.~13, no.~R51, 2006.

\bibitem{firro2007longest}
{\sc Firro, G.}, {\sc Mansour, T.} and {\sc Wilson, M.C.}, ``{Longest alternating subsequences in pattern-restricted permutation},'' {\em the electronic journal of combinatorics}, vol.~14, no.~R34, 2007.

\bibitem{hoeffding1948central}
{\sc Hoeffding, W.} and {\sc Robbins, H.}, ``{The central limit theorem for
  dependent random variables},'' {\em Duke Math. J}, vol.~15, no.~3,
  pp.~773--780, 1948.
  
\bibitem{heinrich1985non}
{\sc Heinrich, L.}, ``{Non-uniform estimates, moderate and large deviations in the central limit theorem for m-dependent random variables},'' {\em Mathematische Nachrichten}, vol.~121, no.~1, pp.~107--121, 1985.

\bibitem{houdré2008longest}
{\sc Houdr\'e, C.} and {\sc Litherland, T.L.}, ``{On the longest increasing subsequence for finite and countable alphabets},'' {\em 
High Dimensional Probability V: The Luminy Volume  IMS Collections 5}, pp.~185--212, 2009.

\bibitem{trevis2}
{\sc Houdr\'e, C.} and {\sc Litherland, T.L.}, ``{On the Limiting Shape of Random Young Tableaux For Markovian Words},'' {\em 
arXiv:0810.2982}, 2009.


\bibitem{lin1996limit}
{\sc Lin, Z.}, {\sc Zhengyan, L.}, {\sc Lu, C.}, and {\sc Chuanrong, L.}, {\em
  {Limit theory for mixing dependent random variables}}.
\newblock Kluwer Academic Pub, 1996.


\bibitem{mansour2008longest}
{\sc Mansour, T.}, ``{Longest alternating subsequences of k-ary words},'' {\em
  Discrete Applied Mathematics}, vol.~156, no.~1, pp.~119--124, 2008.

\bibitem{mansour2008blongest}
{\sc Mansour, T.}, ``{Longest alternating subsequences in pattern-restricted k-ary words},'' {\em Online J. Analytic Combin}, vol.~3, 2008.

\bibitem{resnick1999probability}
{\sc Resnick, S.}, {\em {A probability path}}.
\newblock Birkhauser, 1999.

\bibitem{riauba1977local}
{\sc Riauba, B.},  ``{A local limit theorem for dependent random variables},'' {\em
  Lithuanian Mathematical Journal}, vol.~17, no.~1, pp.~119--129, 1977.

\bibitem{shiryaev1996probability}
{\sc Shiryaev, A.}, ``{Probability. Number 95 in Graduate Texts in
  Mathematics},'' 1996.

\bibitem{stanley511419longest}
{\sc Stanley, R.}, ``{Longest alternating subsequences of permutations},'' {\em
  Michigan Mathematical Journal}, vol.~57, pp.~675--687, 2008.

\bibitem{stanley2006increasing}
{\sc Stanley, R.}, ``{Increasing and decreasing subsequences and their
  variants},'' {\em Proc. Internat. Cong. Math (Madrid 2006)}, American Mathematical Society, pp.~549--579, 2007.

\bibitem{volkonskii1959some}
{\sc Volkonskii, V.} and {\sc Rozanov, Y.}, ``{Some limit theorems for random
  functions. I},'' {\em Theory of Probability and its Applications}, vol.~4,
  p.~178, 1959.

\bibitem{widom2006limiting}
{\sc Widom, H.}, ``{On the Limiting Distribution for the Length of the Longest
  Alternating Sequence in a Random Permutation},'' {\em The Electronic Journal
  of Combinatorics}, vol.~13, no.~R25, 2006.

\bibitem{zeilberger1990holonomic}
{\sc Zeilberger, D}, ``{A holonomic systems approach to special functions identities* 1},'' {\em Journal of computational and applied mathematics}, vol.~32, no.~3, p.~321-368, 1990.

\end{thebibliography}
\end{document}